\magnification=\magstep1 

\def\[{[\![}
\def\]{]\!]}
\def\Z{{\bf Z}}
\def\C{{\bf C}}
\def\Ch{{\C[[h]]}}
\def\l{\ldots}
\def\n{\noindent}

\def\n{\noindent}
\def\t{\theta}

\def\k{{\bar k}}
\def\L{{\bar L}}
\def\q{{\bar q}}

\def\d{\delta}
\def\v{\varepsilon}
\def\ra{\rangle}
\def\la{\langle}

\def\c{\centerline}

\hfuzz=10pt 
\font\tfont=cmbx10 scaled\magstep2 
\font\ttfont=cmbx10 scaled\magstep3 
\font\sfont=cmbx10 scaled\magstep1 
\font\male=cmr9
\vsize=8.5truein \hsize=6truein
\nopagenumbers 
\topskip=0.5truein
\raggedbottom 
\abovedisplayskip=3mm 
\belowdisplayskip=3mm 
\abovedisplayshortskip=0mm 
\belowdisplayshortskip=2mm 
\normalbaselineskip=12pt 
\normalbaselines


\null 
\vskip 1truecm 
 
\centerline{{\ttfont  A  new description of the quantum superalgebra}}

\centerline{{\ttfont U$_{\bf q}$[sl(n+1$|$m)] and related Fock
representations}}

\vskip 1cm

\centerline{{\sfont T.D. Palev$^{\bf a)}$, N.I. Stoilova
$^{\bf b)}$\footnote*{\rm Permanent address: Institute for Nuclear 
Research and Nuclear Energy, Boul. Tsarigradsko Chaussee 72,
1784 Sofia, Bulgaria.} and J. Van der Jeugt$^{\bf c)}$\footnote{**}{\rm 
Research Associate of the Fund for Scientific Research - Flanders (Belgium);
(e-mail: Joris.VanderJeugt@rug.ac.be).} 
}} 

\vskip 0.5 cm 
 
\centerline{$^{\bf a)}$Institute for Nuclear Research and Nuclear 
Energy }
\centerline{Boul. Tsarigradsko Chaussee 72}
\centerline{1784 Sofia, Bulgaria}

\centerline{$^{\bf b)}$Arnold Sommerfeld Institute for Mathematical 
Physics }
\centerline{Technical University of Clausthal}
\centerline{Leibnizstra\ss e 10}
\centerline{D-38678 Clausthal-Zellerfeld, Germany}

\centerline{$^{\bf c)}$Department of Applied Mathematics and Computer
Science }
\centerline{University of Ghent}
\centerline{Krijgslaan 281-S9}
\centerline{B-9000 Gent, Belgium}

\vskip 1cm 

\centerline{\bf Abstract} 

\midinsert\narrower\narrower\male
\n A description of the quantum superalgebra
$U_q[sl(n+1|m)]$ via creation and annihilation generators
(CAGs) is given. A statement that the Fock
representations of the CAGs provide microscopic realizations
of exclusion statistics is formulated.
\endinsert

\vskip 1cm 

\noindent 
{\tfont 1 \ Description of U$_{\bf q}$[sl(n+1$|$m)] via deformed 
CAGs} 
\vskip 0.3cm

\noindent 
First we introduce $U_q[sl(n+1|m)]$ by means of its classical 
definition in terms of the Cartan matrix $\alpha_{ij}$ and the 
Chevalley
generators $h_i, e_i,f_i, \;i,j=1,\ldots ,n+m$.  
Let $\Ch$ be the complex algebra of formal power series 
in the indeterminate $h$, $q=e^h\in \Ch$.
$U_q[sl(n+1|m)]$ is a Hopf algebra, which
is a topologically free $\Ch$ module with generators 
$h_i, e_i,f_i,$ subject to the 
Cartan-Kac relations (${\bar q}=q^{-1}$)
$$
\eqalignno{
& [h_i,h_j]=0,& (1a)\cr
& [h_i,e_j]=\alpha_{ij}e_j,\quad [h_i,f_j]=-\alpha_{ij}f_j,& (1b)\cr
& \[e_i,f_j\]=\delta _{ij}{k_i-\bar{k}_i\over{q-\bar{q}}},\quad
  k_i=q^{h_i},\;k_i^{-1}\equiv\k_i=q^{-h_i}    ,&(1c)\cr 
& \alpha_{ij}=(1+(-1)^{\t_{i-1,i}})\delta_{ij}-
(-1)^{\t_{i-1,i}}\delta_{i,j-1}-\delta_{i-1,j},&\cr
& \t_i=\cases {{ 0}, & if $\; i=0,1,2, \ldots , n$,\cr 
               {1}, & if $\; i=n+1,n+2,\ldots ,n+m$,\cr }; \quad
\t_{ij}=\t_i+\t_j,&\cr
}
$$
the $e$-Serre relations 
$$
\eqalignno{
& [e_i,e_j]=0,\; if \; |i-j|\neq 1;\quad  e^2_{n+1}=0, & (2a) \cr
& [e_i, [e_{i}, e_{i\pm 1}]_{\bar{q}}]_q=
  [e_i, [e_{i}, e_{i\pm 1}]_{q}]_{\bar{q}}=0, \quad i\neq n+1, & 
(2b)\cr
& \{ e_{n+1},[[e_n,e_{n+1}]_q, e_{n+2}]_{\bar{q}}\}=&\cr
&   \{ e_{n+1},[[e_n,e_{n+1}]_{\bar{q}}, e_{n+2}]_{q}\}=0, & (2c)\cr
}
$$
and the $f-$Serre relations, obtained from the $e$-Serre relations
by replacing everywhere
$e_i$ with $f_i$.

\noindent
Here and everywhere: 

\c{$[a,b]_x=ab-xba,~~~  \{a,b\}_x=ab+xba$,}
\c{$\[a,b\]_x=ab-(-1)^{deg(a)deg(b)}xba. $}

We do not write the other Hopf algebra maps $(\Delta,\; \varepsilon,
S)$, since we will not use them. They are certainly also a
part of the definition.

Introduce a normal order in the system of the positive roots [1,2]

\c{$\Delta_+=\{\v_i-\v_j|i<j=0,\l ,n+m\}$}

\n
as follows:
$$\v_i-\v_j<\v_k-\v_l\; \;\;if \;j<l\;\; or \; if\;\; j=l\;\; 
and \;\;i<k. \eqno(3) $$

\n Define the deformed CAGs to be
Cartan-Weyl basis vectors, which are in agreement 
with the above normal order: 
$$
\eqalignno{
& a_1^-=e_1,\quad a_1^+=f_1, \quad H_1=h_1, & \cr
& a_i^-=[[[\l[[e_1,e_2]_{\bar{q}_1},e_3]_{\bar{q}_2},\l]_{\bar{q}_{i-
3}},
e_{i-1}]_{\bar{q}_{i-2}},e_i]_{\bar{q}_{i-1}}=
[a_{i-1}^-, e_i]_{\bar{q}_{i-1}},& \cr
&
a_i^+=[f_i,[f_{i-1},[\l ,[f_3,[f_2,f_1]_{q_1}]_{q_2}\l]_{q_{i-3}}
]_{q_{i-2}}]_{q_{i-1}}=[f_i,a_{i-1}^+]_{q_{i-1}},&\cr
&  H_i=h_1+(-1)^{\t_1}h_2+(-1)^{\t_2}h_3+\l 
+(-1)^{\t_{i-1}}h_i,& (4)\cr
}
$$
where
$$
q_i=q^{1-2\t_i}=\cases{ q,  & if $i\le n$,\cr
                       \q, & if $i>n$.\cr } \eqno(5)
$$
{\it Theorem.} $U_q[sl(n+1|m)]$ is an unital associative
algebra with generators \hfill\break
$H_i,\; a_i^\pm, \;\; i=1,\l ,n+m$
and relations
$$
\eqalignno{
& [H_i,H_j]=0, & (6a) \cr
& [H_i,a_j^{\pm}]=\mp(1+(-1)^{\t_i}\delta_{ij})a_j^{\pm}, 
&(6b)\cr
& \[ a_i^-, a_i^+\]={L_i-\bar{L}_i\over{q-\bar{q}}},\quad 
L_i=q^{H_i},\;\bar{L}_i\equiv L_i^{-1}=q^{-H_i}  & (6c)\cr
& \[\[a_i^{\eta}, a_{i+\xi}^{-\eta}\], a_k^{\eta} 
\]_{q^{\xi (1+(-1)^{\t_i}\delta_{ik})}}=
\eta^{\t_k}\delta_{k,i+\xi}L_k^{-\xi\eta}
a_i^{\eta},  & (6d)\cr
& \[ a_1^\xi , a_2^\xi \]_q =0, \quad \[ a_1^\xi , a_1^\xi \] =0,
\quad \xi, \eta =\pm \;\;or \;\; \pm 1. & (6e) \cr
}
$$

\vskip 1cm 

\noindent 
{\tfont 2 \ Fock representations} 
\vskip 0.3cm

\noindent 
The irreducible Fock representations are
labelled by one non-negative integer $p=1,2,\ldots$, 
called an order of the statistics. To construct them assume that the
corresponding representation space $W_p$ contains 
 a cyclic vector $|0\ra$, such that
$$
a_i^-|0\ra =0,\quad H_i|0\ra =p|0\ra, \quad i=1,2,\ldots ,n+m;
$$
$$
\[ a_{ i}^-,a_{ j}^+\]|0\ra=0,\quad
i\neq j=1,2, \ldots ,n+m. \eqno(7)
$$

From (6) one derives
 that the deformed creation 
(resp. annihilation)
generators $q$-supercommute,
$$
 \[ a_i^\xi , a_j^\xi \]_{q} =0,\quad 
\; i<j=1,\l ,n+m, \;\; \xi=\pm.\eqno(8)
$$
This makes evident the basis (or at least one possible basis) in a 
given Fock space, since any product of only creation generators 
can be always ordered.

As a basis in the Fock space $W_p$ take the vectors
$$
\eqalignno{
& 
|p;r_1,r_2,\ldots,r_{n+m})=
 \sqrt{[p-\sum_{l=1}^{n+m} r_l]!
\over {[p]![r_1]!\ldots[r_{n+m}]!}}
(a_1^+)^{r_1}(a_2^+)^{r_2}\ldots (a_n^+)^{r_n} \times  & \cr
& (a_{n+1}^+)^{r_{n+1}}
(a_{n+2}^+)^{r_{n+2}}
\ldots (a_{n+m}^+)^{r_{n+m}}|0\ra,~~[x]={q^x-q^{-x}\over{q-q^{-1}}},
\quad q\in {\bf R} & (9)\cr
}
$$
with
$$
r_i\in \Z_+,\;\; i=1,\l ,n; \quad r_{i}\in \{0,1\}, \;\;
i=n+1,\l ,n+m,~~~ 
\sum_{i=1}^{n+m}r_i\leq p. 
\eqno(10)
$$
In order to write down the transformations of the basis under the
action of the CAG's one has to write down the supercommutation
relations between all Cartan-Weyl generators, expressed via
the CAGs. The latter is a necessary condition for the application of
the Poincare-Birghoff-Witt theorem, when computing the action of
the generators on the Fock basis vectors. 
Bellow we display only some necessary relations, which follow
from (6) in a rather nontrivial way:
$$
\eqalignno{
&  L_i\L_i=\L_iL_i=1, \quad  L_iL_j=L_jL_i, \quad
L_ia_j^\pm=q^{\mp (1+ (-1)^{\t _i}\delta_{ij})}a_j^\pm L_i, & 
(11)\cr
& \[a_i^-,a_i^+\]={L_i-{\bar L}_i\over q-{\bar q}},\quad
\[a_i^\eta,a_j^\eta\]_q=0,\quad \eta=\pm,\quad i<j, & (12)  \cr 
}
$$
$$
\eqalignno{
& \[\[a_i^\eta,a_j^{-\eta}\],a_k^\eta\]_{q^{\xi(1+(-
1)^{\t_i}\d_{ik}})}=\eta^{\t_j}\d_{jk}L_k^{-\xi\eta}a_i^\eta + 
(-1)^{\t_k}\epsilon(j,k,i)
(q-\q)\[a_k^\eta,a_j^{-\eta}\]a_i^\eta=&  \cr
& \eta^{\t_j}\d_{jk}L_k^{-\xi\eta}a_i^\eta + 
(-1)^{\t_k\t_j}\epsilon(j,k,i)q^\xi
(q-\q)a_i^\eta\[a_k^\eta,a_j^{-\eta}\], \; \xi(j-i)>0, \;
\xi,\;\eta =\pm & (13)\cr
}
$$
where
$$
\epsilon(j,k,i)=\cases {\;\;\;1, & if $j>k>i$;\cr -1, & if 
$j<k<i$;\cr 
\;\;\; 0, & otherwise.\cr}
\eqno (14)
$$
{\it Proposition.} The set of all vectors (9)
constitute an orthonormal basis in  $W_p$ with respect
to the scalar product, defined in the usual way with "bra" and
"ket" vectors and $\la 0|0 \ra =1$. 
The transformation of the basis under the
action of the CAOs read: 
$$
\eqalignno{
& H_i|p;r_1,r_2,\ldots,r_{n+m})=\left(p-(-1)^{\t_i}r_i-
\sum_{j=1}^{n+m}r_j
\right) |p;r_1,r_2,\ldots,r_{n+m}), & (15)\cr
& a_i^-|p;r_1,\ldots,r_{n+m})=
(-1)^{\t_i(\t_1r_1+\l +\t_{i-1}r_{i-1})}q^{r_1+\ldots +r_{i-
1}}\sqrt{[r_i]
[p-\sum_{l=1}^{n+m} r_l +1]}  &\cr
& \times |p;r_1,\ldots r_{i-1},r_i-1,r_{i+1},
\ldots,r_{n+m}),& (16)   \cr
&  a_i^+|p;r_1,\ldots,r_{n+m})=
(-1)^{\t_i(\t_1r_1+\l +\t_{i-1}r_{i-1})}
\q^{r_1+\ldots +r_{i-1}}(1-\t_ir_i)
\sqrt{[r_i+1][p-\sum_{l=1}^{n+m} r_l]}&\cr
& \times |p;r_1,\ldots r_{i-1},r_i+1,r_{i+1},\ldots,r_{n+m}).& (17)   
\cr
}
$$

\vskip 1cm 

\noindent 
{\tfont 3 \  Properties of the underlying statistics} 
\vskip 0.3cm

\noindent 
In the present talk we have defined the algebra 
$U_q[sl(n+1|m)]$ in terms of new set of generators, called creation
and annihilation generators. Let us illustrate on a simple
example of the nondeformed algebra $sl(n+1|m)$ that within each
Fock representation  $a_i^+$ (resp. $a_i^-$) can be interpreted
as operators creating (resp. annihilating) "particles"
with, say, energy $\varepsilon_i$. Only for simplicity let
us assume that $n=m$. Set moreover
$b_i^\pm=a_i^\pm,~~ 
f_i^\pm=a_{i+n}^\pm,~~ i=1,\ldots,n,$
and consider a "free"
Hamiltonian 
$
H=\sum_{i=1}^{n}\varepsilon_i (H_i + H_{i+n})=
\sum_{i=1}^{n}\varepsilon_i (\[b_i^+,b_i^-\]+\[f_i^+,f_i^-\]).$
Then 
$[H,b_i^\pm]=\pm \varepsilon_i b_i^\pm,~~~
[H,f_i^\pm]=\pm \varepsilon_i f_i^\pm. $
This result together with (nondeformed) Eqs. (16) and (17) allows one 
to interprete $r_i, ~i=1,\ldots,n$  as the number of $b-$particles 
with 
energy
$\varepsilon_i$ and similarly $r_{i+n}, ~i=1,\ldots,n$  as the number 
of 
$f-$particles with energy $\varepsilon_i$.
Then $b_i^+$ ($f_i^+$) increases this number by
one, it ``creates'' a particle in the one-particle state (=
orbital) $i$.  Similarly, the operator $b_i^-$ ($f_i^-$) diminishes 
this
number by one, it ``kills'' a particle on the $i-th$ orbital. On
every orbital $i$ there cannot be more than one particle of kind $f$,
whereas such restriction does not hold for the $b-$particles.
These are, so to speak, Fermi like (resp. Bose like) properties.
There is however one essential difference. If the order of the
statistics is $p$ than no more than $p$ "particles" can be
accommodated in the system,
$\sum_{i=1}^{n+m}r_i\leq p.$
Hence the 
available places for new particles to be "loaded" on, say, $i^{th}$
orbital 
depend on how many particles 
have been already accomodated on the other orbitals.
This is neither Bose, nor Fermi like property. It is however
a typical property for the so called {\it exclusion statistics} [3].

The statistics, which we have outlined  above,  
is indeed an exclusion statistics. 
This statement will be proved rigorously elsewhere. 
The exclusion statistics reformulates the concept of fractional
statistics as a 
{\it generalized Pauli exclusion principle} for
spaces with arbitrary dimensions.  Despite of the fact that
exclusion statistics are defined for any space dimensions, so far
they were applied and tested only within 1D and 2D models. The
statistics described above are examples of microscopic
description of exclusion statistics in arbitrary space
dimensions.

\vskip 1cm 

\noindent 
{\tfont Acknowledgments}
\vskip 0.3cm

\n
T.D.P. is grateful to Prof. H.D. Doebner for the invitation
to attend the International Symposium on Quantum Theory and
Symmetries, 18-22 July 1999, Goslar. N.I.S. is thankful to
Prof. H.D. Doebner for the kind hospitality to work at ASI, TU
Clausthal and to DAAD for the three months fellowship.

\vskip 1cm 

\noindent
{\tfont References} 
\vskip 0.3cm

\n
1. V.N. Tolstoy, {\it Lecture Notes in Physics}, {\bf 370}, 
   (Berlin: Springer) p. 118 (1990).

\n
2. S.M. Khoroshkin and V.N. Tolstoy, {\it Commun. Math. Phys.},
   {\bf 141}, 599 (1991). 

\n
3. F.D.M. Haldane, {\it Phys. Rev. Lett.}, {\bf 67},  937 (1991).

\end